\documentclass[12pt]{article}

\usepackage{latexsym}
\usepackage{amsfonts}
\usepackage{amssymb}
\usepackage{amsmath}
\usepackage{mathrsfs}
\usepackage{bbm}
\usepackage{dsfont}
\usepackage{hyperref}

\newcommand{\be}{\begin{equation}}
\newcommand{\ee}{\end{equation}}
\newcommand{\bea}{\begin{eqnarray}}
\newcommand{\eea}{\end{eqnarray}}
\newcommand{\bean}{\begin{eqnarray*}}
\newcommand{\eean}{\end{eqnarray*}}
\newcommand{\brray}{\begin{array}}
\newcommand{\erray}{\end{array}}

\newtheorem{dfn}{Definition}[section]
\newtheorem{thm}[dfn]{Theorem}
\newtheorem{lmma}[dfn]{Lemma}
\newtheorem{ppsn}[dfn]{Proposition}
\newtheorem{crlre}[dfn]{Corollary}
\newtheorem{xmpl}[dfn]{Example}
\newtheorem{rmrk}[dfn]{Remark}

\newcommand{\bdfn}{\begin{dfn}\rm}
\newcommand{\bthm}{\begin{thm}}
\newcommand{\blmma}{\begin{lmma}}
\newcommand{\bppsn}{\begin{ppsn}}
\newcommand{\bcrlre}{\begin{crlre}}
\newcommand{\bxmpl}{\begin{xmpl}}
\newcommand{\brmrk}{\begin{rmrk}\rm}

\newcommand{\edfn}{\end{dfn}}
\newcommand{\ethm}{\end{thm}}
\newcommand{\elmma}{\end{lmma}}
\newcommand{\eppsn}{\end{ppsn}}
\newcommand{\ecrlre}{\end{crlre}}
\newcommand{\exmpl}{\end{xmpl}}
\newcommand{\ermrk}{\end{rmrk}}

\newcommand{\clh}{\mathcal{H}}

\newcommand{\clk}{\mathcal{K}}

\author{ S. Sundar}
\title { An asymmetric multiparameter CCR flow}

\begin{document}
\maketitle
\begin{abstract}
 In this  note, we exhibit an example of a multiparameter CCR flow which is not cocycle conjugate to its opposite. This is in sharp contrast to the
 one parameter situation. 
 \end{abstract}
\noindent {\bf AMS Classification No. :} {Primary 46L55; Secondary 46L99.}  

\section{Opposite of an $E_0$-semigroup}
Let $P \subset \mathbb{R}^{d}$ be a closed convex cone. We assume that $d \geq 2$, $P$ spans $\mathbb{R}^{d}$ and $P$ contains no line, i.e. $P \cap -P=\{0\}$. We denote the interior of $P$ by $\Omega$.  For $x,y \in \mathbb{R}^{d}$, we write $x \leq y$ ($x<y$) if $y-x \in P$ ($y-x \in \Omega$). Let $\alpha:=\{\alpha_{x}\}_{x \in P}$ be an $E_0$-semigroup over $P$ on $B(\clh)$ where $\clh$ is an infinite dimensional separable Hilbert space. For $x \in P$, let 
\[
E(x):=\{T \in B(\clh): \alpha_{x}(A)T=TA ~~\textrm{for all $A \in B(\clh)$}\}.
\]
For $x \in P$, $E(x)$ is a separable Hilbert space where the inner product is given by $\langle T|S \rangle=S^{*}T$. The disjoint union  of Hilbert spaces $\displaystyle \mathcal{E}:=\coprod_{x \in \Omega}E(x)$ has a structure of a product system and is called the product system associated with $\alpha$. It is indeed a cocycle conjugacy invariant. For more on product systems and $E_0$-semigroups in the multiparameter context, we refer the reader to \cite{Murugan_Sundar_continuous}.

Keep the foregoing notation. Let \[
\mathcal{E}^{op}:=\{(x,T) \in \Omega \times B(\clh): x \in \Omega, T \in E(x)\}.\] 
Define a semigroup multiplication on $\mathcal{E}^{op}$ by the following formula:
\begin{equation}
\label{product for opposite}
(x,T).(y,S)=(x+y,ST)\end{equation}
for $(x,T), (y,S) \in \mathcal{E}^{op}$. Then $\mathcal{E}^{op}$ is an abstract product system over $\Omega$ (in the sense of Definition 2.1 of \cite{Murugan_Sundar_continuous}). Arveson's bijection between product systems and $E_0$-semigroups, established in \cite{Murugan_Sundar_continuous} for the case of a cone, ensures that there exists an $E_0$-semigroup denoted $\alpha^{op}:=\{\alpha^{op}_{x}\}_{x \in P}$, which is unique up to cocycle conjugacy, such that the product system associated to $\alpha^{op}$ is isomorphic to $\mathcal{E}^{op}$. The $E_0$-semigroup $\alpha^{op}$ is  called the opposite of $\alpha$.

A natural question that arises in this context is the following.  Are the $E_0$-semigroups $\alpha$ and $\alpha^{op}$ cocycle conjugate ?
If $\alpha$ is cocycle conjugate to $\alpha^{op}$, we call $\alpha$ symmetric otherwise we call $\alpha$ asymmetric.
In the one parameter context, we have the following.
\begin{enumerate}
\item[(1)] One parameter CCR flows are symmetric. This follows from Arveson's classification of type I $E_0$-semigroups and the index computation.
\item[(2)] Tsirelson in his remarkable papers \cite{Tsi} and \cite{Tsirelson2003} constructed examples of type $II_0$ $E_0$-semigroups which are asymmetric by probabilistic means.
\end{enumerate}
It is natural to ask whether $(1)$ stays true in the multiparameter context. More precisely, suppose  $d \geq 2$ and $V$ is a pure isometric representation of $P$ on a Hilbert space $\clh$. Let $\alpha^{V}$ be the CCR flow associated to $V$ acting on $B(\Gamma(\clh))$ where $\Gamma(\clh)$ is
the symmetric Fock space of $\clh$. Is $\alpha^{V}$ symmetric ? We show that for the left regular representation of $P$ on $L^{2}(P)$, the associated CCR flow is asymmetric.

\section{Decomposable product systems}
Following Arveson, the author in \cite{Injectivity} has defined the notion of a decomposable product system. Let us review the definitions. Let $\alpha:=\{\alpha_{x}\}_{x \in P}$ be an $E_0$-semigroup on $B(\clh)$ and let
$E:=\{E(x)\}_{x \in P}$ be the associated product system. Suppose $x \in P$ and $u \in E(x)$ is a non-zero vector. We say that $u$ is decomposable if given $y \leq x$ with $y \in P$, there exists $v \in E(y)$ and $w \in E(x-y)$ such
that $u=vw$. The set of decomposable vectors is denoted by $D(x)$. We say that $\alpha$ is decomposable if 
\begin{enumerate}
\item[(1)] for $x,y \in P$, $D(x)D(y)=D(x+y)$, and
\item[(2)] for $x \in P$, $D(x)$ is total in $E(x)$.
\end{enumerate}

\begin{ppsn}
\label{the opposite}
The opposite of a decomposable $E_0$-semigroup is decomposable.
\end{ppsn}

\begin{rmrk}
The subtle point that we wish to stress is that apriori the product rule in the opposite product system given by Eq. \ref{product for  opposite} holds only over $\Omega$. However to prove Prop. \ref{the opposite}, we need the validity of the product rule over the whole semigroup $P$ which is assured by the next lemma. 

\end{rmrk}

Let us fix a few notation. Let $\alpha:=\{\alpha_{x}\}_{x \in P}$ be a decomposable $E_0$-semigroup acting on $B(\clh)$. The product system of $\alpha$ is denoted by $E:=\{E(x)\}_{x \in P}$. We denote the opposite of $\alpha$ by $\beta:=\{\beta_{x}\}_{x \in P}$. 
Suppose that $\beta$ acts on $B(\widetilde{\clh})$. Denote the product system of $\beta$ by $F:=\{F(x)\}_{x \in P}$.

\begin{lmma}
\label{antiunitary}
For every $x \in P$, there exists a unitary  $\widetilde{\theta}_{x}:E(x) \to F(x)$ such that for $x,y \in P$, $T \in E(x)$ and $S \in E(y)$, 
\[
\widetilde{\theta}_{x+y}(ST)=\widetilde{\theta}_{x}(T)\widetilde{\theta}_{y}(S).\]
\end{lmma}
\textit{Proof.} From the definition of $\beta$, it follows that for $x \in \Omega$, there exists a unitary operator $\theta_{x}:E(x) \to F(x)$ such that
for $x,y \in \Omega$, $T \in E(x)$ and $S \in E(y)$, 
\[
\theta_{x+y}(ST)=\theta_{x}(T)\theta_{y}(S).\]
 Fix an element $a \in \Omega$. Let $x \in P$ and $T \in E(x)$ be given.
We claim that there exists a unique bounded linear operator on $\widetilde{\clh}$, which we denote by $\widetilde{\theta_{x}}(T)$, such that
for $S \in E(a)$ and $\xi \in \widetilde{\clh}$,
\[
\widetilde{\theta}_{x}(T)\theta_{a}(S)\xi=\theta_{a+x}(ST)\xi.\]

 For any $b \in \Omega$, the map $E(b) \otimes \widetilde{\clh} \ni S \otimes \xi \to \theta_{b}(S)\xi \in \widetilde{\clh}$ is a unitary operator. This way, we can identify $\widetilde{\clh}$ with $E(b) \otimes \widetilde{\clh}$ for every $b \in \Omega$. 
Right multiplication by $T$ induces a bounded linear operator from $E(a) \to E(a+x)$ of norm $||T||$. Tensor with identity to obtain the desired operator $\widetilde{\theta}_{x}(T)$ from $\widetilde{\clh} \equiv E(a) \otimes \widetilde{\clh} \to E(a+x) \otimes \widetilde{\clh}\equiv \widetilde{\clh}$. This proves our claim.

As the set $\{\theta_{a}(S)\xi: S \in E(a), \xi \in \widetilde{\clh}\}$ is total in $\widetilde{\clh}$ and $\{\theta_{x}\}_{x \in \Omega}$ is anti multiplicative, it follows that $\widetilde{\theta}_{x}=\theta_{x}$ if $x \in \Omega$. Let $x \in P$, $y \in \Omega$, $T \in E(x)$ and $S \in E(y)$ be given. For $R \in E(a)$ and $\xi \in \widetilde{\clh}$, calculate as follows to observe that
\begin{align*}
\widetilde{\theta}_{x+y}(TS)\theta_{a}(R)\xi&=\theta_{(a+x)+y}(RTS)\xi\\
                                                                    &=\theta_{y}(S)\theta_{a+x}(RT)\xi \\
                                                                    &=\theta_{y}(S)\widetilde{\theta_x}(T)\theta_{a}(R)\xi.
 \end{align*}
Hence it follows that 
\begin{equation}
\label{Eq 1}
\widetilde{\theta}_{x+y}(TS)=\theta_{y}(S)\widetilde{\theta_x}(T)
\end{equation} for $x \in P$, $y \in \Omega$, $T \in E(x)$ and $S \in E(y)$. 

Let $x,y \in P$, $T \in E(x)$ and $S \in E(y)$ be given.
 Let $R \in E(a)$ be of unit norm. Calculate as follows to observe 
that
\begin{align*}
\theta_{a}(R)\widetilde{\theta}_{x+y}(TS)&=\theta_{x+y+a}(TSR) ~~(\textrm{by Eq. \ref{Eq 1}}) \\
                                                                 &= \theta_{x+(y+a)}(T(SR)) \\
                                                                 &=\theta_{y+a}(SR)\widetilde{\theta}_{x}(T)~~(\textrm{by Eq. \ref{Eq 1}}) \\
                                                                 &=\theta_{a}(R)\widetilde{\theta}_{y}(S)\widetilde{\theta}_{x}(T)~~(\textrm{by Eq. \ref{Eq 1}}) .
\end{align*}
Premultiplying the above equation by $\theta_{a}(R)^{*}$, we get
\begin{equation}
\label{Eq 2}
\widetilde{\theta}_{x+y}(TS)=\widetilde{\theta}_{y}(S)\widetilde{\theta}_{x}(T)
\end{equation}
for $x,y \in P$, $T \in E(x)$ and $S \in E(y)$. 

Note that for $x \in P$, the map $E(x) \ni T \to \widetilde{\theta}_{x}(T) \in B(\widetilde{\clh})$ is linear and norm preserving. A direct calculation reveals that for $x \in P$, $T_1,T_2 \in E(x)$, 
\begin{equation}
\label{Eq 3}
\widetilde{\theta}_{x}(T_2)^{*}\widetilde{\theta}_{x}(T_1)=\langle T_1|T_2 \rangle_{E(x)}.
\end{equation}
For $x \in P$, let $\widetilde{F}(x):=\{\widetilde{\theta}_{x}(T): T \in E(x)\}$.  Fix $x \in P$. It follows from  Eq. \ref{Eq 3}  that there exists a unique normal $*$-endomorphism denoted $\widetilde{\beta}_{x}$ on $B(\widetilde{\clh})$ such that the intertwining space of $\widetilde{\beta}_{x}$ is $\widetilde{F}(x)$. Eq. \ref{Eq 2} implies that the family $\widetilde{\beta}:=\{\widetilde{\beta}_{x}\}_{x \in P}$ forms a semigroup of endomorphisms. 

Note that $\widetilde{F}(x)=F(x)$ for every $x \in \Omega$. Hence $\widetilde{\beta}_{x}=\beta_{x}$ for $x \in \Omega$. The semigroup $\widetilde{\beta}$ agrees with an $E_0$-semigroup $\beta$ on the interior $\Omega$. Thanks to Lemma 4.1 of \cite{Murugan_Sundar}, it follows that $\widetilde{\beta}$ is an $E_0$-semigroup. Since $\Omega$ is dense in $P$, it follows that $\widetilde{\beta}_{x}=\beta_{x}$ for every $x \in P$. Consequently, we have $\widetilde{F}(x)=F(x)$ for every $x \in P$. This completes the proof. \hfill $\Box$

\textbf{Proof of Proposition \ref{the opposite}}: Let $\{\widetilde{\theta}_{x}\}_{x \in P}$ be a family of unitary operators as in Lemma \ref{antiunitary}. For $x \in P$, denote the set of decomposable vectors in $E(x)$ by $D(x)$. A moment's reflection on the definition shows  that the decomposable vectors of $F(x)$ is $\{\widetilde{\theta}_{x}(T): T \in D(x)\}$. The conclusion is now immediate. \hfill $\Box$

\section{ A counterexample}
In this section, we produce the promised counterexample, i.e. a CCR flow which is asymmetric. Let us recall the definition of a CCR flow associated to a pure isometric representation. Suppose $V:P \to B(\clh)$ is a pure isometric representation. Recall that $V$ is said to be pure if  $\bigcap_{x \in P}V_{x}\clh=\{0\}$. Denote the symmetric Fock space of $\clh$ by $\Gamma(\clh)$. The CCR flow associated to $V$, denoted $\alpha^{V}:=\{\alpha_{x}\}_{x \in P}$, is the unique $E_0$-semigroup
on $B(\Gamma(\clh))$ such that the following equation is satisfied. For $x \in P$ and $\xi \in \clh$,
\[
\alpha_{x}(W(\xi))=W(V_{x}\xi)
\]where $\{W(\xi):\xi \in \clh\}$ is the set of Weyl operators on $\Gamma(\clh)$. For more details regarding multiparameter CCR flows, we refer the reader to \cite{Anbu} and \cite{Anbu_Sundar}.

\begin{rmrk}
\label{isometry}
A few remarks are in order.
\begin{enumerate}
\item[(1)] In \cite{Injectivity}, a strongly continuous isometric representation, indexed by $\Omega$, is constructed out of a decomposable $E_0$-semigroup (see Proposition 4.1 of \cite{Injectivity}). Moreover the resulting isometric representation, up to unitary equivalence, is a cocycle conjugacy invariant. 
\item[(2)] It is shown in \cite{Injectivity} that the CCR flow $\alpha^{V}$ is decomposable. Moreover the isometric representation constructed out of the decomposable $E_0$-semigroup $\alpha^{V}$ is $V$ itself (see Proposition 5.1 of \cite{Injectivity}). 
\end{enumerate}
\end{rmrk}
Fix a pure isometric representation $V$ of $P$ on a Hilbert space $\clh$. Denote the CCR flow $\alpha^{V}$ by $\alpha$ and its opposite by $\alpha^{op}$. Denote the isometric representation constructed out of $\alpha^{op}$ by $V^{op}$.  If we unwrap all the details regarding the construction of $V^{op}$, which is routine, we see that $V^{op}$ has the following explicit description.  Denote the Hilbert space on which $V^{op}$ acts by $\clh^{op}$.

Define an equivalence relation on  $\{(\xi,a): \xi \in Ker(V_{a}^{*}), a \in \Omega\}$ as follows. We say $(\xi,a) \sim (\eta,b)$ if and only if $V_{b}\xi=V_{a}\eta$. Let $H^{op}$ be the set of equivalence classes. Then $H^{op}$ has an inner product space structure where the addition, scalar multiplication and inner product are given by
\begin{align*}
[(\xi,a)]+[(\eta,b)]&=[(V_{b}\xi+V_{a}\eta,a+b)]\\
\lambda[(\xi,a)]&=[(\lambda\xi,a)] \\
\langle [(\xi,a)]|[(\eta,b)]\rangle&=\langle V_{b}\xi|V_{a}\eta \rangle.
\end{align*}
Then $\clh^{op}$ is the completion of $H^{op}$ and for $a \in \Omega$, the operator $V_{a}^{op}$ is given by the equation
\[
V_{a}^{op}[(\xi,b)]=[(\xi,a+b)].\]

To produce a counter example of a CCR flow which is asymmetric, it suffices to construct a pure isometric representation $V$ such that $V$ and $V^{op}$ are not unitarily equivalent. For if $\alpha:=\alpha^{V}$ and $\alpha^{op}$ are cocycle conjugate then by Remark \ref{isometry}, it follows that $V$ and $V^{op}$ are unitarily equivalent.

First we obtain a better description of  $V^{op}$. Let $U$ be the minimal unitary dilation of $V$. More precisely, there exists a Hilbert space $\widetilde{\clh}$ containing $\clh$ as a closed subspace and a strongly continuous unitary representation
$U:=\{U_{x}\}_{x \in \mathbb{R}^{d}}$ on $\widetilde{\clh}$ such that the following conditions are satisfied.
\begin{enumerate}
\item[(1)] For $a \in \Omega$ and $\xi \in \clh$, $U_{a}\xi=V_{a}\xi$, and
\item[(2)] the increasing union $\displaystyle \bigcup_{a \in \Omega}U_{a}^{*}\clh$ is dense in $\widetilde{\clh}$.
\end{enumerate}
The minimal unitary dilation is unique up to unitary equivalence and the existence of such a dilation is given by an inductive limit procedure. 

Set $\clk:=\clh^{\perp}$. Since $\clh$ is invariant under $\{U_{a}\}_{a \in \Omega}$, it follows that $\clk$ is invariant under
$\{U_{-a}:a \in \Omega\}=\{U_{a}^{*}: a \in \Omega\}$. For $a \in \Omega$, let $W_{a}$ be the operator on $\clk$ which is the restriction of $U_{-a}$. Then $W:=\{W_{a}\}_{a \in \Omega}$ is a strongly continuous semigroup of isometries 
on $\clk$.  

\begin{ppsn}
\label{Opposite of isometric representation}
With the foregoing notation, we have the following. 

\begin{enumerate}
\item[(1)] The isometric representation $W$ is pure, i.e. $\displaystyle \bigcap_{a \in \Omega}W_{a}\clk=\{0\}$.
\item[(2)] The isometric representations $W$ and $V^{op}$ are unitarily equivalent. 
\end{enumerate}
\end{ppsn}
\textit{Proof.} Fix a point $a \in \Omega$. Set $S:=V_{a}$. Recall the following Archimedean property.  Given $x \in \mathbb{R}^{d}$, there exists a positive integer $n$ such that $na>x$ (see Lemma 3.1 of \cite{Murugan_Sundar_continuous}). 
Thus $\displaystyle \{0\}=\bigcap_{b \in \Omega}V_{b}\clh=\bigcap_{n \geq 1}V_{na}\clh=\bigcap_{n \geq 1}S^{n}\clh.$
In other words, $S$ is a pure isometry. 

By the Archimedean prinicple, we have the equaltiy $\displaystyle \bigcup_{b \in \Omega}U_{b}^{*}\clh=\bigcup_{n \geq 1}U_{a}^{*n}\clh$. Hence $\displaystyle \bigcup_{n \geq 1}U_{a}^{*n}\clh$ is dense in $\widetilde{\clh}$. In otherwords, the discrete one parameter group of unitaries $\{U_{na}: n \in \mathbb{Z}\}$ is the minimal unitary dilation of the discrete one parameter isometric representation $\{S^{n}\}_{ n \geq 0}$. 

Using Wold decomposition, we can identify $\clh$ with $\ell^{2}(\mathbb{N}) \otimes K$ for some Hilbert space $K$ and $S$ with the standard one sided shift with multiplicity.  Then $\widetilde{\clh}$ can be identified with $\ell^{2}(\mathbb{Z}) \otimes K$ with $U_{a}$ identified with the bilateral shift with multiplicity. Once this identification is made, it is clear that $\displaystyle \bigcap_{n \geq 0}W_{na}\clk=\bigcap_{n \geq 0}U_{a}^{*n}\clk=\{0\}$. Once again by the Archimedean principle, we have the equality
$
\displaystyle \bigcap_{b \in \Omega}W_{b}\clk=\bigcap_{n \geq 0}W_{na}\clk=\{0\}$.
This proves $(1)$. 

Let $a \in \Omega$. We claim that the image of the map
\[
Ker(W_{a}^{*}) \ni \xi \to U_{a}\xi \in \widetilde{\clh}\]  is contained in  $Ker(V_{a}^{*})$. Let $\xi \in Ker(W_{a}^{*})$ and $\eta \in \clk$ be given.  Observe that 
\[
\langle U_{a}\xi|\eta \rangle =\langle \xi|U_{a}^{*}\eta \rangle 
                                            =\langle \xi|W_{a}\eta \rangle 
                                            =\langle W_{a}^{*}\xi|\eta \rangle 
                                            =0.
\]
This proves that for $\xi \in Ker(W_{a}^{*})$, $U_{a}\xi \in \clh$. Let $\xi \in Ker(W_{a}^{*})$ and $\eta \in \clh$ be given. Using the fact that $\xi \in \clk$ and $\eta \in \clh$, observe that
\[
\langle U_{a}\xi|V_{a}\eta\rangle=\langle U_{a}\xi|U_{a}\eta\rangle=\langle \xi|\eta \rangle=0.\]
Thus $U_{a}\xi$ is orthogonal to $Ran(V_{a})$. Coupled with the fact that $U_{a}\xi \in \clh$, we conclude that $U_{a}\xi \in Ker(V_{a}^{*})$. 
This proves our claim. 
A calculation similar to the one above implies that the image of the map 
$
Ker(V_{a}^{*}) \ni \eta \to U_{-a}\eta \in \widetilde{\clh}$ is contained in $Ker(W_{a}^{*})$. Consequently the map
$
Ker(W_{a}^{*}) \ni \xi \to U_{a}\xi \in Ker(V_{a}^{*})$
is a unitary. 

By $(1)$, we have $\clk=\overline{\bigcup_{a \in \Omega}Ker(W_{a}^{*})}$. Note that the family of inner product preserving maps \[\Big\{Ker(W_{a}^{*}) \ni \xi \to [(U_{a}\xi,a)] \in \clh^{op}\Big\}_{a \in \Omega}\] patch together and defines a unitary map from $\clk=\overline{\bigcup_{a \in \Omega}Ker(W_{a}^{*})}$ to $\clh^{op}$, which we denote by $T$, such that the following holds. 
If $\xi \in Ker(W_{a}^{*})$ for some $a \in \Omega$ then \[
 T\xi=[(U_{a}\xi,a)].\]
 It is clear that $T$ intertwines the isometric representations $W$ and $V^{op}$. This proves $(2)$. The proof is now complete. \hfill $\Box$

Let $A \subset \mathbb{R}^{d}$ be such that $A$ is non-empty, proper, closed and  $A+P \subset A$. Such subsets were called $P$-modules in \cite{Anbu_Sundar}. Consider the Hilbert space $L^{2}(A)$. For $x \in P$, let $V_{x}$ be the isometry on $L^{2}(A)$
defined by the equation
\begin{equation}
\label{isometries}
V_{x}(f)(y):=\begin{cases}
 f(y-x)  & \mbox{ if
} y -x \in A,\cr
   &\cr
    0 &  \mbox{ if } y-x \notin A.
         \end{cases}
\end{equation}
Then $(V_{x})_{x \in P}$ is an isometric representation of $P$ which we denote by $V^{A}$. We call $V^{A}$ the isometric representation associated to the $P$-module $A$. Moreover the representation $V^{A}$ is pure. In what follows, $Int(A)$ denotes the interior of $A$.

\begin{lmma}
\label{dilation of shift}
Let $A$ be a $P$-module.  We have the following.
\begin{enumerate}
\item[(1)] The increasing union $\bigcup_{a \in \Omega}(Int(A)-a)=\mathbb{R}^{d}$.
\item[(2)] Given a compact subset $K \subset \mathbb{R}^{d}$, there exists $a \in \Omega$ such that $K$ is contained in $Int(A)-a$.
\item[(3)] The minimal unitary dilation of $V^{A}$ is the left regular representation of $\mathbb{R}^{d}$ on $L^{2}(\mathbb{R}^{d})$. 
\end{enumerate}
\end{lmma}
\textit{Proof.} Since $A$ is a $P$-module, it is clear that if $a<b$ with $a,b \in \Omega$ then $Int(A)-a $ is contained in $Int(A)-b$. By translating, if necessary, we can assume $0 \in A$. Hence $P \subset A$ and $\Omega \subset Int(A)$. Observe the equality
$\displaystyle \mathbb{R}^{d}=\Omega-\Omega=\bigcup_{a \in \Omega}(\Omega-a) \subset \bigcup_{a \in \Omega}(Int(A)-a)$. This proves $(1)$. 

Fix an interior point $a \in \Omega$. By $(1)$ and by the Archimedean property, it follows that $(Int(A)-na)_{n \geq 1}$ is an increasing sequence of open sets which increases to $\mathbb{R}^{d}$. Now $(2)$ is immediate.

Let $\{U_{x}\}_{x \in \mathbb{R}^{d}}$ be the left regular representation of $\mathbb{R}^{d}$ on $L^{2}(\mathbb{R}^{d})$.  From the definition it follows that for $a \in \Omega$, $V_{a}$ is the compression of $U_{a}$ onto $L^{2}(A)$. Observe  that for $a \in \Omega$, $U_{a}^{*}L^{2}(A) =L^{2}(A-a)$. It follows from $(2)$ that $\displaystyle \bigcup_{a \in \Omega}U_{a}^{*}L^{2}(A)$ contains the space of continuous functions with compact support. Consequently $\displaystyle \bigcup_{a \in \Omega}U_{a}^{*}L^{2}(A)$ is dense in $L^{2}(\mathbb{R}^{d})$. Hence $\{U_{x}\}_{x \in \mathbb{R}^{d}}$ has all the properties required to be the minimal  minimal unitary dilation of $V^{A}$. This completes the proof. \hfill $\Box$

Note that if $A$ is a $P$-module then $-(Int A)^{c}$ is a $P$-module. 
Fix a $P$-module $A$  and let $V:=V^{A}$ be the isometric representation associated to $A$. Set $B:=-(Int(A))^{c}$.
\begin{ppsn}
 Keep the foregoing notation. 
 \begin{enumerate}
 \item[(1)]  The isometric representation $V^{op}$ is unitarily equivalent to the representation $V^{B}$.
\item[(2)] The representations $V$ and $V^{op}$ are unitarily equivalent if and only if there exists $z \in \mathbb{R}^{d}$ such that $A=B+z$. 
\end{enumerate}
\end{ppsn}
\textit{Proof.} Thanks to Lemma II.12 of \cite{Hilgert_Neeb}, $A$ and $Int(A)$ differ by a set of measure zero. Thus we can identify $L^{2}(A)$ with $L^{2}(Int(A))$. Lemma \ref{dilation of shift} and Proposition \ref{Opposite of isometric representation} implies that $V^{op}$ is equivalent to the isometric representation $W=\{W_{a}\}_{a \in \Omega}$ acting on the Hilbert space $L^{2}((Int(A))^{c})$ where the operators $\{W_{a}\}_{a \in \Omega}$ are given by the following equation. For $a \in \Omega$ and $f \in L^{2}((Int(A))^{c})$
\begin{equation}
\label{isometries}
W_{a}(f)(y):=\begin{cases}
 f(y+a)  & \mbox{ if
} y +a \in (Int(A))^{c},\cr
   &\cr
    0 &  \mbox{ if } y+a \notin (Int(A))^{c}.
         \end{cases}
\end{equation}
The inversion $\mathbb{R}^{d} \ni x \to -x \in \mathbb{R}^{d}$ induces a unitary between the Hilbert spaces $L^{2}((Int(A))^{c})$ and $L^{2}(B)$ and intertwines the representations $W$ and $V^{B}$. This proves $(1)$.

It is clear that if $A$ is a translate of $B$ then $V=V^{A}$ and $V^{op}=V^{B}$ are unitarily equivalent. On the other hand, suppose $V^{A}$ and $V^{B}$ are unitarily equivalent. Then the associated CCR flows $\alpha^{V^{A}}$ and $\alpha^{V^{B}}$ are 
cocycle conjugate. By Theorem 1.2 of \cite{Anbu_Sundar}, it follows that $A$ and $B$ are translates of each other.\footnote{See also Page 26 of \cite{Injectivity}.} This completes the proof. \hfill $\Box$

Thus to produce a counterexample of a CCR flow which is not cocycle conjugate to its opposite, it suffices to produce a $P$-module $A$ such that $A$ is not a translate of $-(Int(A))^{c}$. The cone $P$ itself is one such candidate. Recall that we have assumed $d \geq 2$. Let us recall the notion of extreme points.
For a subset $C $ of $\mathbb{R}^{d}$ and a point $x \in C$, we say $x$ is an extreme point of $C$ if  $x=\frac{y+z}{2}$ with $y,z \in C$ then $y=z=x$.  

\begin{lmma}
The sets $P$ and $-\Omega^{c}$ are not translates of each other. 
\end{lmma}
\textit{Proof.} First we claim that $P \cup -P\neq \mathbb{R}^{d}$. Suppose not. Since we have assumed that $P \cap -P =\{0\}$, it follows that  the only boundary point of $P$ is $\{0\}$. Fix $a \in \Omega$. Proposition 2.3 of \cite{Anbu_Sundar} implies that the map
\[
\partial(P) \times (0,\infty) \ni (x,s) \to x+sa \in \Omega\]
is a homeomorphism. But  $\partial(P)=\{0\}$ and hence  $\Omega=\{sa: s>0\}$. Since $\Omega$ spans $\mathbb{R}^{d}$, it follows that $d=1$  contradicting our assumption. Therefore $P \cup -P \neq \mathbb{R}^{d}$.

Note that the set of extreme points of $P$ is $\{0\}$. For we have assumed that $P$ contains no line. On the other hand, we claim that $\Omega^{c}$ has no extreme point. Note that $t\Omega^{c} \subset \Omega^{c}$ for $t>0$. Hence the set of extreme points of $\Omega^{c}$ is contained in $\{0\}$. But $0$ is not an extreme point of $\Omega^{c}$. To see this, pick $x \notin P \cup -P$. Then $x \in \Omega^{c}$, $-x \in \Omega^{c}$ and $x \neq 0$. But $0=\frac{x+(-x)}{2}$. This proves that the set of extreme points of $\Omega^{c}$ is empty. 
Consequently, the set of extreme points of any translate of $-\Omega^{c}$ is empty. Hence $P$ and $-\Omega^{c}$ are not translates of each other. This completes the proof. \hfill $\Box$.

\begin{rmrk}
 Let $V$ be a pure isometric representation of $P$ and $\alpha=\alpha^{V}$ be the associated CCR flow. Let $V^{op}$ be the isometric representation corresponding to the decomposable $E_0$-semigroup $\alpha^{op}$. Since $\alpha$ is spatial, it follows that  $\alpha^{op}$ is spatial. (Recall that an $E_0$-semigroup is said to be spatial if  its product system admits a nowhere vanishing multiplicative measurable cross section).
 By Theorem 4.4 of \cite{Injectivity}, it follows that $\alpha^{op}$ and $\alpha^{V^{op}}$ are cocycle conjugate. Thus, to summarise, an opposite of a CCR flow is a CCR flow but not necessarily the same as the original one. 
 
 \end{rmrk}
 
\nocite{Arveson}
\nocite{Tsi}
\nocite{Tsirelson2003}

\bibliography{references}
 \bibliographystyle{amsplain}

\noindent
{\sc S. Sundar}
(\texttt{sundarsobers@gmail.com})\\
         {\footnotesize  Institute of Mathematical Sciences (HBNI), CIT Campus, \\
Taramani, Chennai, 600113, Tamilnadu, INDIA.}\\

\end{document}